\documentclass{amsproc}

\usepackage{graphicx} 

\usepackage{}

\newtheorem{theorem}{Theorem}[section]

\theoremstyle{definition}
\newtheorem{definition}[theorem]{Definition}
\newtheorem{example}[theorem]{Example}

\newtheorem{problem}[theorem]{Problem}

\newtheorem{conjecture}[theorem]{Conjecture}

\theoremstyle{remark}
\newtheorem{remark}[theorem]{Remark}
\newtheorem{question}[theorem]{Question}

\numberwithin{equation}{section}

\begin{document}

\title{Multibranched surfaces in 3-manifolds}


\author{Makoto Ozawa}
\address{Department of Natural Sciences, Faculty of Arts and Sciences, Komazawa University, 1-23-1 Komazawa, Setagaya-ku, Tokyo, 154-8525, Japan}
\curraddr{}
\email{w3c@komazawa-u.ac.jp}
\thanks{The author is partially supported by Grant-in-Aid for Scientific Research (C) (No. 17K05262) and (B) (No. 16H03928), The Ministry of Education, Culture, Sports, Science and Technology, Japan}



\date{}

\begin{abstract}
This is a latest survey article on embeddings of multibranched surfaces into 3-manifolds.
\end{abstract}

\maketitle


Throughout this article, we will work in the piecewise linear category.
All topological spaces are assumed to be second countable and Hausdorff.
It is a fundamental problem that for two topological spaces $X$ and $Y$, 
\begin{enumerate}
\item Can $X$ be embedded in $Y$?
\item If $X$ can be embedded in $Y$, then
\begin{enumerate}
\item When are two embeddings of $X$ in $Y$ equivalent?
\item How can $X$ be embedded in $Y$?
\end{enumerate}
\end{enumerate}
In this article, we consider the case that $X$ is a multibranched surface and $Y$ is a closed orientable 3-manifold.

We say that a 2-dimensional CW complex is a {\em multibranched surface} if
we remove all points whose open neighborhoods are homeomorphic to the
2-dimensional Euclidean space, then we obtain a 1-dimensional complex
which is homeomorphic to a disjoint union of some simple closed
curves.

Multibranched surfaces naturally arise in several areas:
\begin{itemize}
\item Polycontinuous patterns|a mathematical model of
microphase-separated structures made by block copolymers (\cite{HCO}, \cite{HS}, \cite{FCKHS})
\item 2-stratifolds | as spines of closed 3-manifolds (\cite{GGH}, \cite{GGH16}, \cite{GGH17})
\item Trisections, multisections | as an analogue of Heegaard splittings (\cite{GK}, \cite{Ko}, \cite{RT})
\item Essential surfaces | as non-meridional essential surfaces in link exteriors (\cite{EM}, \cite{EO}), essential surfaces in handlebody-knot exteriors (\cite{KO2}) and Dehn surgery (\cite{EM}, \cite{GL})
\end{itemize}

In Section 1, we define several terms related with multibranched
surfaces.
In Section 2, we describe several backgrounds of multibranched surfaces.
In Section 3, we study embeddings of multibranched surfaces into closed orientable 3-manifolds.
In Section 4, we consider on multibranched surfaces which cannot be
embedded into the 3-sphere.

\section{Definitions}

\subsection{Definition}

\begin{figure}[htbp]
	\begin{center}
	\includegraphics[width=0.4\linewidth]{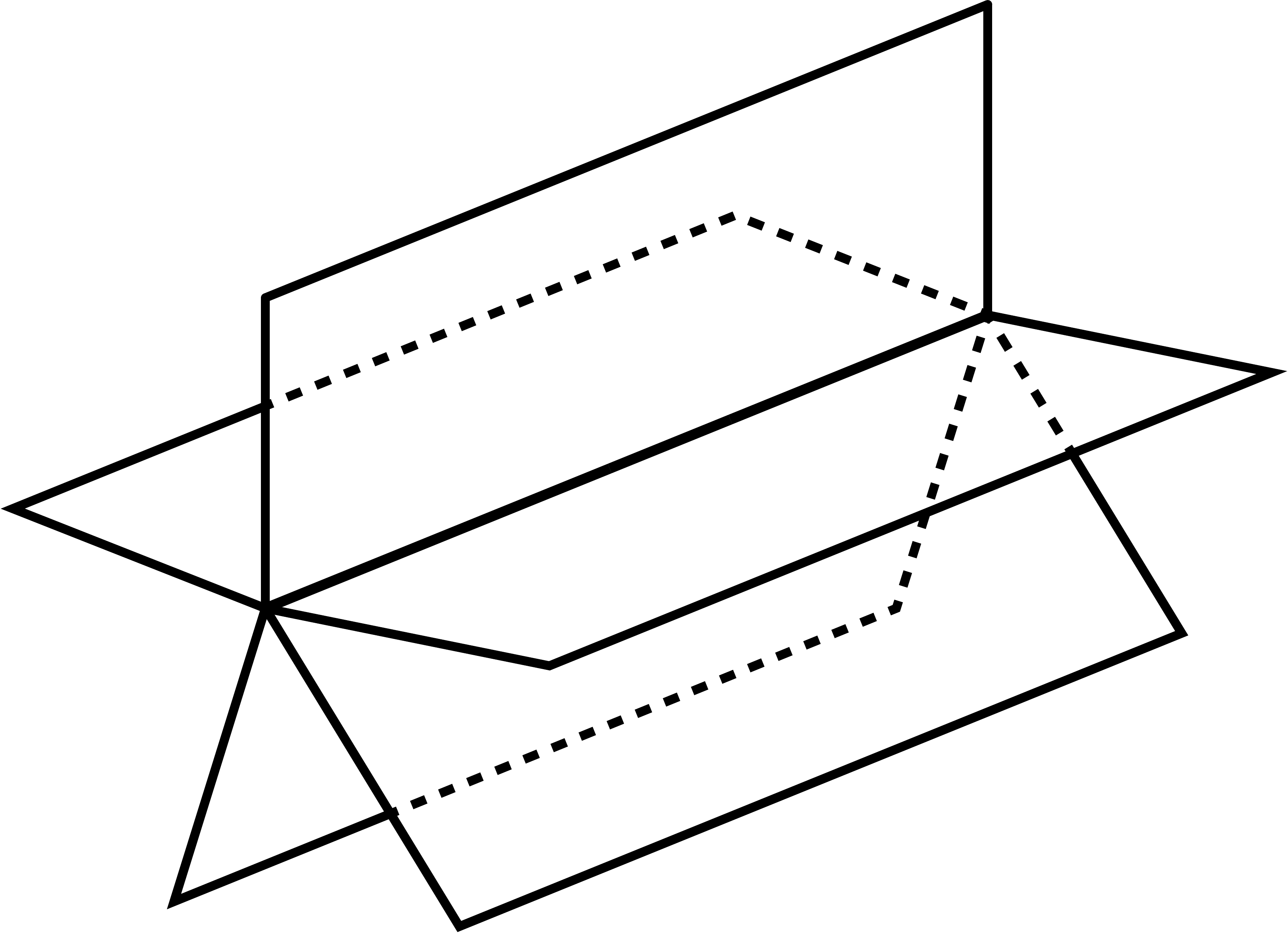}
	\end{center}
	\caption{The multibranched Euclidean plane $\Bbb{R}^2_5$}
	\label{model}
\end{figure}

Let $\Bbb{R}^2_+$ be the closed upper half-plane $\{(x_1,x_2)\in \Bbb{R}^2 \mid x_2\ge 0\}$.
The {\em multibranched Euclidean plane}, denoted by $\Bbb{R}^2_i$ $(i\ge 1)$, is the quotient space obtained from $i$ copies of $\Bbb{R}^2_+$ by identifying with their boundaries $\partial \Bbb{R}^2_+=\{(x_1,x_2)\in\Bbb{R}^2\mid x_2=0\}$ via the identity map.
See Figure \ref{model} for the multibranched Euclidean plane $\Bbb{R}^2_5$.

A second countable Hausdorff space $X$ is called a {\em multibranched surface} if $X$ contains a disjoint union of simple closed curves $l_1,\ldots, l_n$ satisfying the following:
\begin{enumerate}
\item For each point $x\in l_1\cup \cdots \cup l_n$, there exist an open neighborhood $U$ of $x$ and a positive integer $i$ such that $U$ is homeomorphic to $\Bbb{R}^2_i$.
\item For each point $x\in X-(l_1\cup\cdots\cup l_n)$, there exists an open neighborhood $U$ of $x$ such that $U$ is homeomorphic to $\Bbb{R}^2$.
\end{enumerate}

\subsection{Construction}\label{construction}

To construct a compact multibranched surface, 
we prepare a closed $1$-dimensional manifold $B$ (corresponding to $l_1,\ldots, l_n$), a compact $2$-dimensional manifold $S$ (corresponding to the closure of each component of $X-(l_1\cup\cdots\cup l_n)$) and a map $\phi: \partial S \to B$ satisfying that
for every connected component $c$ of $\partial S$, the restriction $ \phi|_c : c \to \phi(c) $ is a covering map.
Then a multibranched surface $X$ can be constructed from the triple $(B, S; \phi)$ as the quotient space $X=B \cup_\phi S$.

A connected component of $B$, $S$ or $\partial S$ is said to be a {\em branch},  {\em sector} or {\em prebranch} respectively. 
The set consisting of all branches or sectors is denoted by $\mathcal{B}(X)$ or $\mathcal{S}(X)$ respectively.

\subsection{Degree, oriented degree, regularity}

For a prebranch $c$ of a multibranched surface $X$, the covering degree of $\phi|_c:c\to \phi(c)$ is called the {\em degree} of $c$ and denoted by $d(c)$. 
We give an orientation for each branch and each prebranch $c$ of $X$ (In the case that every sector $s$ is orientable, the orientation of $c$ is induced by that of $s$). 
The {\em oriented degree} of a prebranch $c$ of $X$ is defined as follows: if the covering map $\phi|_c:c\to \phi(c)$ is orientation preserving, the {\em oriented degree} $od(c)$ of $c$ is defined by $od(c)=d(c)$ and if it is orientation reversing, the oriented degree is defined by $od(c)=-d(c)$.

A prebranch $c$ of $X$ is said to be {\it attached} to a branch $l$ if $\phi(c)=l$.
We denote by $\mathcal{A}(l)$ the set consisting of all prebranches which are attached to a branch $l$ and the number of elements of $\mathcal{A}(l)$ is called the {\em index} of $l$ and denoted by $i(l)$.

A multibranched surface $X$ is {\em regular} if for every branch $l$ and every prebranch $c,c'\in\mathcal{A}(l)$, $d(c)=d(c')$.
Let $X$ be a regular multibranched surface. Since each pair of prebranches $c,c'\in\mathcal{A}(l)$ has same degree, the {\em degree} of a branch $l$ is well-defined as $d(l)=d(c)=d(c')$.

\subsection{Graph representation}

Let $X$ be a compact multibranched surface obtained from $(B, S; \phi)$ such that all components of $S$ are orientable and have non-empty boundary (Hereafter, we assume such conditions on multibranched surfaces unless otherwise stated).
The multibranched surface $X=B \cup_\phi S$ has a graph representation (\cite{EO}) as follows.
Let $G=(V_S\cup V_B,E)$ be a bipertite graph such that $|V_S|=|\mathcal{S}(X)|$ and $|V_B|=|\mathcal{B}(X)|$.
For each sector $s\in\mathcal{S}(X)$, we assign a vertex $v(s)\in V_S$ with a label $g(s)$.
For each branch $l\in\mathcal{B}(X)$, we assign a vertex $v(l)\in V_B$.
For a prebranch $c\subset \partial s$, we assign an edge $e\in E$ with a label $od(c)$ connecting $v(s)$ and $v(l)$ where $c\in\mathcal{A}(l)$.
Almost similar concept to this graph representation has been defined in \cite{GGH16}.

\begin{example}\label{crosscap}
A closed non-orientable surface of crosscap number $h$ can be regarded as a multibranched surface $X$ with $h$ branches $B=l_1\cup \cdots \cup l_h$ and a planar surface $S$ with $h$ boundary components, such that $od(c)=2$ for any prebranch $c\subset \partial S$.
Then $X$ has a graph representation $G$ as shown in Figure \ref{non-orientable}.

\begin{figure}[htbp]
	\begin{center}
	\begin{tabular}{cc}
	\includegraphics[width=.3\linewidth]{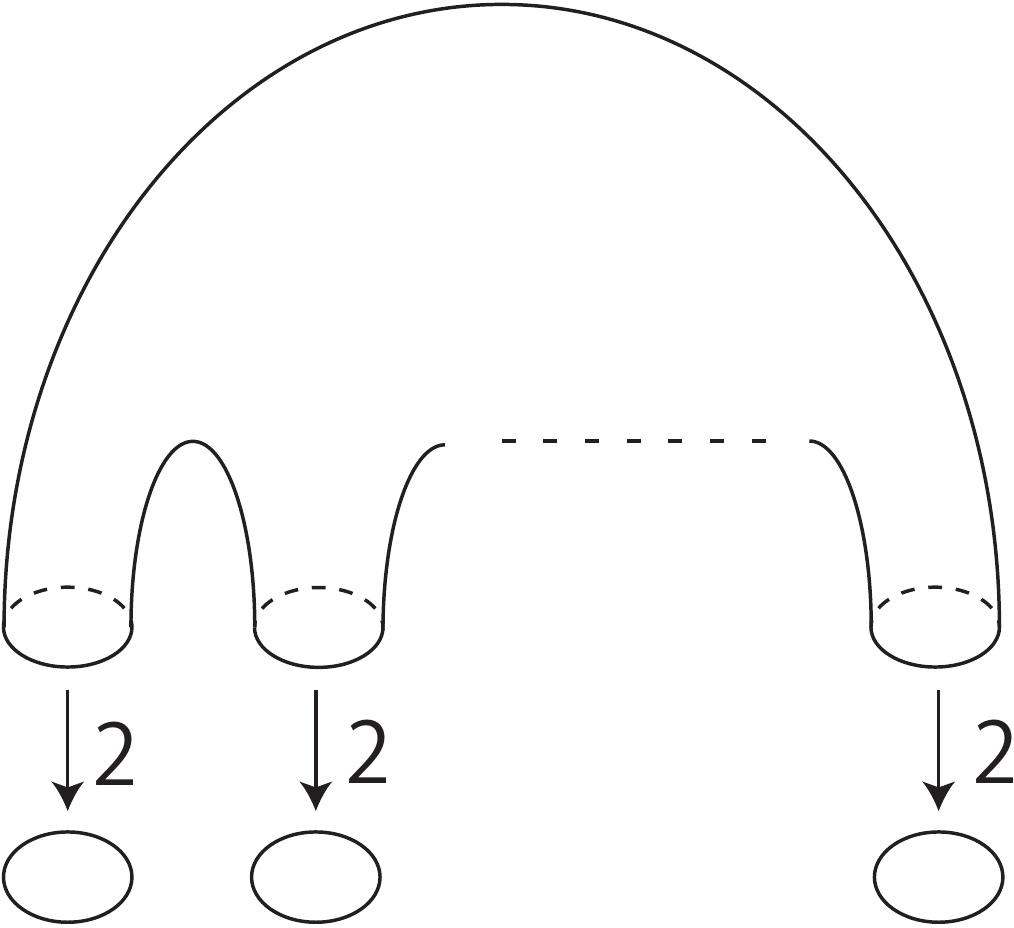}&
	\includegraphics[width=.27\linewidth]{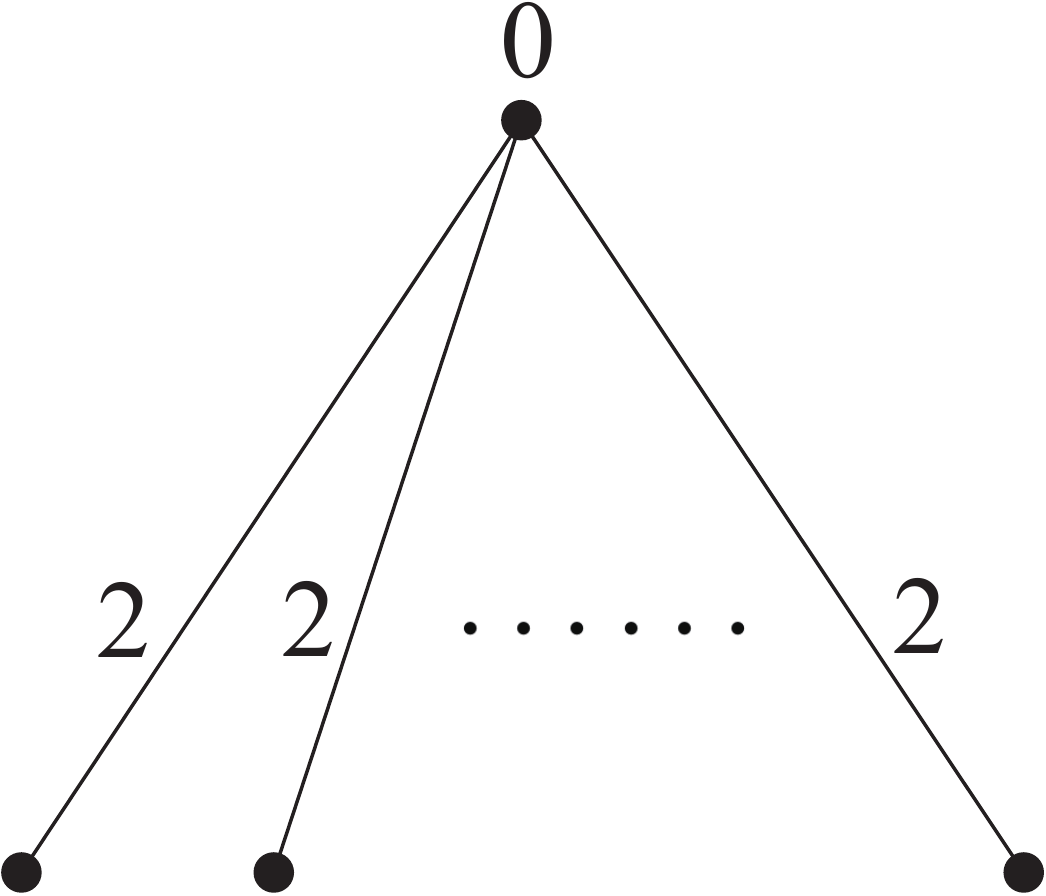}\\
	$X$ & $G$
	\end{tabular}
	\end{center}
	\caption{A multibranched surface $X$ and its graph representation $G$}
	\label{non-orientable}
\end{figure}
\end{example}

\subsection{Incident matrix}

For a sector $s\in\mathcal{S}(X)$ and a branch $l\in\mathcal{B}(X)$, we define the {\em algebraic degree} $ad_s(l)$ as follows.
\[
\displaystyle ad_s(l) = \sum_{c\in \mathcal{A}(l) \cap \partial s} od(c)
\]
Then we define the {\em incident matrix} $M_X=(a_{ij})$ $(i=1,\ldots,n; j=1,\ldots,m)$ by
\[
a_{ij} = ad_{s_j} (l_i),
\]
where $\mathcal{B}(X)=\{l_1,\ldots,l_n\}$ and $\mathcal{S}(X)=\{s_1,\ldots,s_m\}$.

\subsection{First homology group}

For a branch $l$ and a sector $s$ of a regular multibranched surface $X$, we define $d(l;s)=\sum_{c \subset \partial s} od(c)$, where $c$ is a prebranch attached to $l$.
The multibranched surface obtained by the removing an open disk from each sector is denoted by $\dot{X}$.

\begin{theorem}[{\cite[Theorem 4.1]{MO}}]\label{homology}
Let $X$ be a regular multibranched surface with $\mathcal{B}(X)=\left\{ l_1, \ldots, l_n \right\}$, $\mathcal{S}(X)=\left\{ s_1, \ldots, s_m \right\}$.
Then, 
\[
H_1(X) = \left[ l_1, \ldots, l_n : \sum_{k=1}^{n}d(l_k;s_1)l_k, \ldots, \sum_{k=1}^{n}d(l_k;s_m)l_k \right] \oplus \mathbb{Z}^{r'(X)} 
\] 
where $r'(X)=rank H_1(\dot{X})-n$.
\end{theorem}

Therefore, the torsion subgroup of $H_1(X)$ can be calculated from the incident matrix $M_X$.

\begin{example}
Let $X$ be a multibranched surface which has a graph representation as shown in Figure \ref{X_2}, where $n=4$, $g_i=0$ and all wrapping number is 1.
In \cite[Example 4.2]{MO}, the first homology group is calculated by using Theorem \ref{homology} as $H_1(X)=(\Bbb{Z}/3\Bbb{Z})\oplus \Bbb{Z}^4$.

As we shall see later, the incident matrix of $X$ is:
\[
M_{X}=\left( \begin{array}{cccc}
0 & 1 & 1 & 1  \\
1 & 0 & 1 & 1  \\
1 & 1 & 0 & 1  \\
1 & 1 & 1 & 0 \\
\end{array} \right)
\]
This matrix is equivalent to $(3)$ as follows.
\[
\left( \begin{array}{cccc}
0 & 1 & 1 & 1  \\
1 & 0 & 1 & 1  \\
1 & 1 & 0 & 1  \\
1 & 1 & 1 & 0 \\
\end{array} \right)
\sim
\left( \begin{array}{cccc}
3 & 3 & 3 & 3  \\
1 & 0 & 1 & 1  \\
1 & 1 & 0 & 1  \\
1 & 1 & 1 & 0 \\
\end{array} \right)
\sim
\left( \begin{array}{cccc}
3 & 0 & 0 & 0  \\
1 & -1 & 0 & 0  \\
1 & 0 & -1 & 0  \\
1 & 0 & 0 & -1 \\
\end{array} \right)
\]
\[
\sim
\left( \begin{array}{cccc}
3 & 0 & 0 & 0  \\
0 & -1 & 0 & 0  \\
0 & 0 & -1 & 0  \\
0 & 0 & 0 & -1 \\
\end{array} \right)
\sim
\left( \begin{array}{cccc}
3 & 0 & 0 & 0  \\
0 & 1 & 0 & 0  \\
0 & 0 & 1 & 0  \\
0 & 0 & 0 & 1 \\
\end{array} \right)
\sim
\left( \begin{array}{c}
3  \\
\end{array} \right)
\]
This shows that the torsion subgroup of $H_1(X)$ is $\Bbb{Z}/3\Bbb{Z}$.
\end{example}

On the other hand, a natural presentation for the fundamental group of a 2-stratifold was given in \cite{GGH17}.
Thus, we can also obtain the first homology group via abelianization.

\subsection{Circular permutation system and slope system}

A {\em permutation} of a set $A$ is a bijection from the additive group $\mathbb{Z} /n \mathbb{Z}$ into $A$. Two permutations $\sigma$ and $\sigma'$ of $A$ are {\em equivalent} if there is an element $k \in \mathbb{Z} /n \mathbb{Z}$ such that $\sigma'(x) = \sigma (x+k)$ ($x \in \mathbb{Z} /n \mathbb{Z}$). An equivalent class of a permutation of $A$ is a {\em circular permutation}.

For a regular multibranched surface $X$, we define the ``circular permutation system'' and ``slope system'' of $X$ as follows.
A circular permutation of $\mathcal{A}(l)$ is called a {\em circular permutation} on a branch $l$.
A collection $\mathcal{P}=\{ \mathcal{P}_l \}_{l \in \mathcal{L}(X)}$ is called a {\em circular permutation system} of $X$ if $\mathcal{P}_l$ is a circular permutation on $l$.
For a branch $l$, a rational number $p/q$ with $q=d(l)$ is called a {\em slope} of $l$.
A collection $\{ \mathcal{S}_l \}_{l \in \mathcal{L}(X)}$ is called a {\em slope system} of $X$ if $\mathcal{S}_l$ is a slope of $l$.

\subsection{Neighborhoods}


Let $X=B \cup_\phi S$ be a regular multibranched surface and let $\mathcal{P}=\{ \mathcal{P}_l \}_{l \in \mathcal{L}(X)}$ and $\mathcal{S}=\{ \mathcal{S}_l \}_{ l \in \mathcal{L}(X) }$ be a permutation system and a slope system of $X$ respectively. 
We will construct a compact and orientable $3$-manifold, which is uniquely determined up to homeomorphism by a pair of $\mathcal{P}$ and $\mathcal{S}$, by the following procedure. 

First, for each branch $l\in\mathcal{B}(X)$ and each sector $s\in\mathcal{S}(X)$, we take a solid torus $l \times D^2$, where $D^2$ is a disk and take a product $s \times [-1, 1]$. If $s$ is non-orientable, then we take a twisted $I$-bundle $s \tilde{\times} [-1, 1]$ over $s$. We give orientations for these $3$-manifolds. 

Next, we glue them together depending on the permutation system $\mathcal{P}$ and the slope system $\mathcal{S}$, where we assign the slope $\mathcal{S}_l$ of $l$ to the isotopy class of a loop $k_l$ in $\partial ( l \times D^2)$,
by an orientation reversing map $\Phi: \partial S \times [-1, 1] \to \partial (B \times D^2)$ satisfying that for every branch $l$ and every prebranch $c$ with $\phi(c)=l$, the restriction $\Phi |_{c \times [-1, 1]} : c \times [-1, 1] \to N \left(k_l; \partial \left( l \times D^2 \right) \right)$ is a homeomorphism.

Then, we uniquely obtain a compact and orientable $3$-manifold with boundary, denoted by $N(X; \mathcal{P}, \mathcal{S})$. The $3$-manifold $N(X; \mathcal{P}, \mathcal{S})$ is called the {\em neighborhood} of $X$ with respect to $\mathcal{P}$ and $\mathcal{S}$. The set consisting of all neighborhoods of $X$ is denoted by $\mathcal{N}(X)$.

\section{Background}

\subsection{Graphs}

A graph $G$ can be regarded as a 1-dimensional CW complex, where a vertex and an edge correspond to a 0-cell and 1-cell respectively, and vertices of an edge specify the attaching map for a 1-cell to 0-cells.
This structure can be extended to 2-dimensional objects as in Subsection \ref{construction}, that is, we extend vertices, edges, the attaching map to a closed $1$-dimensional manifold $B$ (branch), a compact $2$-dimensional manifold without closed component $S$ (sector), a covering map $\phi: \partial S \to B$ respectively.
Then a multibranched surface $X$ can be obtained as the quotient space $X=B\cup_{\phi}S$.

Kuratowski (\cite{K}) proved that a graph $G$ as a 1-dimensional CW complex cannot be embedded into $\mathbb{R}^2$ if and only if $G$ contains the complete graph $K_5$ or the complete bipartite graph $K_{3,3}$ as a subspace.
At the present time, this result is stated that $G$ cannot be embedded into $\mathbb{R}^2$ if and only if $G$ has $K_5$ or $K_{3,3}$ as a minor.
Robertson--Seymour (\cite{RS}) showed that for any minor closed property $P$, the set of minor minimal graphs which do not have $P$ is finite.
This motivates us to consider the following problem: 
Characterize all ``minor minimal'' multibranched surfaces which cannot be embedded in $\mathbb{R}^3$ (Problem \ref{obstruction_set}).
Since all closed non-orientable surfaces are minor minimal multibranched surfaces, the set of ``minor minimal'' multibranched surfaces which cannot be embedded in $\mathbb{R}^3$ is infinite.
We will describe the details in Section \ref{forbidden}.

\subsection{2-dimensional complexes}

A 2-dimensional CW complex is a multibranched surface if we remove all points whose open neighborhoods are homeomorphic to $\mathbb{R}^2$, then we obtain a 1-dimensional complex which is homeomorphic to a disjoint union of some simple closed curves.
Thus the set of multibranched surfaces is a subset of the set of 2-dimensional CW complexes.

Embeddings of 2-complexes into manifolds are widely studied in \cite{HM}.

Matou\u{s}ek--Sedgwick--Tancer--Wagner (\cite{MSTW}) showed that there is an algorithm that, given a 2-dimensional simplicial complex $K$, decides whether $K$ can be embedded (piecewise linearly, or equivalently, topologically) in $\mathbb{R}^3$.

Carmesin (\cite{C}, \cite{C2}, \cite{C3}, \cite{C4}, \cite{C5}) proved that a locally 3-connected simply connected 2-dimensional simplicial complex has a topological embedding into 3-space if and only if it has no space minor from a finite explicit list $\mathcal{Z}$ of obstructions.

\subsection{Essential surfaces}\label{essential}

The embedding of multibranched surfaces in the 3-sphere $S^3$ is closely related to the existence of essential surfaces in link exteriors.
Let $L$ be a link in $S^3$ and $F$ be an essential surface properly embedded in the exterior $E(L)$ of $L$ whose boundary is non-meridional.
By shrinking the regular neighborhood $N(L)$ into $L$ and extending $F$ along it, we obtain an essential multibranched surface $X$ embedded in $S^3$, where we say that a multibranched surface $X$ with branches $B$ and sectors $S$ embedded in $S^3$ is {\em essential} if $S\cap E(B)$ is essential, namely, incompressible, boundary-incompressible and not boundary-parallel in $E(B)$.
Conversely, let $X$ be an essential multibranched surface with branches $B$ and sectors $S$ embedded in $S^3$.
Then $B$ is a link in $S^3$ and $S\cap E(B)$ is an essential surface properly embedded in $E(B)$ whose boundary is non-meridional.
Therefore, the set of all pairs of a link in $S^3$ and an essential surface properly embedded in the exterior of it whose boundary is non-meridional is equal to the set of all essential multibranched surfaces embedded in $S^3$.

\subsection{Fundamental problem}

The Menger--N\"{o}beling theorem ({\cite[Theorem 1.11.4.]{E}}) shows that any finite $2$-dimensional CW complex can be embedded in $\mathbb{R}^{5}$. 
Furthermore, any multibranched surface can be embedded in $\mathbb{R}^4$ (\cite[Proposition 2.3]{MO}).
More generally, any finite 2-dimensional simplicial complex whose intrinsic 1-skeleton is a proper subset of $K_7$ embeds in $\mathbb{R}^4$ (\cite{G}).

If for a branch $l$, there exist prebranches $c, c'\in \mathcal{A}(l)$ such that $d(c)\ne d(c')$, then the multibranched surface cannot be embedded in any 3-manifold.
The converse also holds, namely we have shown that a multibranched surface can be embedded in some closed orientable 3-manifold if and only if the multibranched surface is regular (\cite[Corollary 2.4]{RBS}, \cite[Proposition 2.7]{MO}).

We remark that any 3-manifold can be embedded in $\Bbb{R}^5$ (\cite{W}).
Thus, we obtain the next diagram on the embedability of multibranched surfaces (Figure \ref{embedding}).

\begin{figure}[htbp]
	\begin{center}
	\includegraphics[width=.85\linewidth]{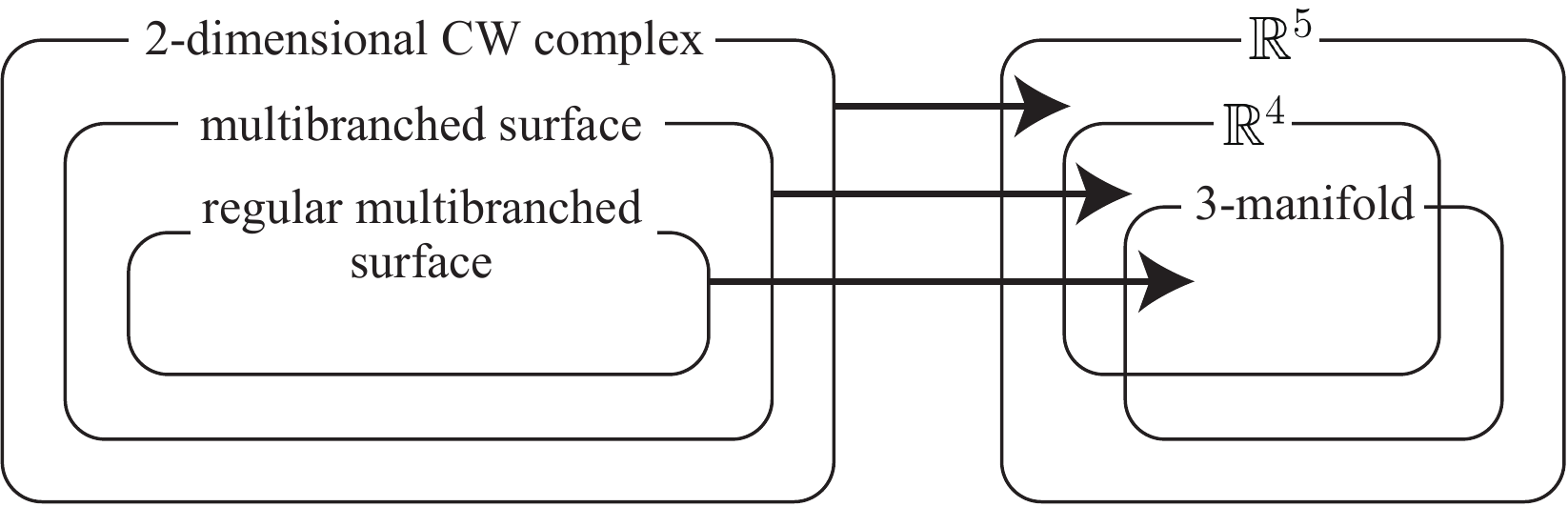}
	\end{center}
	\caption{The embeddability of multibranched surfaces}
	\label{embedding}
\end{figure}

The following problems are fundamental for embeddings of multibranched surfaces.

\begin{problem}\label{3-manifold}
For a regular multibranched surface $X$, find a most simple closed orientable 3-manifold $M$ in which $X$ can be embedded.
Moreover, determine the minimal Heegaard genus of such 3-manifold $M$.
\end{problem}

\begin{problem}\label{3-sphere}
For a regular multibranched surface $X$, determine whether or not $X$ can be embedded in the 3-sphere $S^3$.
\end{problem}

We consider Problem \ref{3-manifold} in Section \ref{embedding}, and Problem \ref{3-sphere} in Section \ref{forbidden}.

\section{Embeddings into 3-manifolds}\label{embedding}

\subsection{Genus}

For a closed orientable 3-manifold $M$, the Heegaard genus is a fundamental index.
The {\em Heegaard genus} $g(M)$ of $M$ is defined as the minimal genus of a closed orientable surface embedded in $M$ which separates $M$ into two orientable handlebodies.

For an orientable compact 3-manifold $N$ with boundary, the minimal Heegaard genus of closed orientable 3-manifolds in which $N$ can be embedded is denoted by $eg(N)$ and called the {\em embeddable genus} of $N$.
We remark that $eg(N)\le g(N)$ (\cite[Proposition 3.1]{MO}), where $g(N)$ denotes the minimal genus of Heegaard splittings of $N$ in a sense of Casson--Gordon (\cite{CG}).

For a regular multibranched surface $X$, we define the {\em minimum genus} $\min g(X)$ and {\em maximum genus} $\max g(X)$ respectively as follows.
\[
\min g(X)=\min \{ eg(N) \mid N \in \mathcal{N}(X) \}
\]
\[
\max g(X)=\max \{ eg(N) \mid N \in \mathcal{N}(X) \}
\]

\subsection{Upper bounds}
There are following inequalities which give upper bounds for the minimum genus and the maximum genus.
The next theorem states originally that $\min g(X) \le |\mathcal{B}(X)| + |\mathcal{S}(X)|$, but its proof is still effective for $\max g(X)$ and implies the latter half.

\begin{theorem}[{\cite[Theorem 3.5]{MO}}]
For a regular multibranched surface $X$, the next inequality holds.
\[
\max g(X) \le |\mathcal{B}(X)| + |\mathcal{S}(X)|
\]
Moreover, if the wrapping number of each branch is 1, then
\[
\max g(X) \le |\mathcal{S}(X)|.
\]
\end{theorem}

\begin{remark}
In the proof of {\cite[Theorem 3.5]{MO}}, it is shown that $X$ can be embedded in a connected sum of $|\mathcal{B}(X)|$ lens spaces and $|\mathcal{S}(X)|$ $S^2\times S^1$'s.
Yuya Koda asked me whether any closed orientable 3-manifold contains a minimal genus
embedding of some multibranched surfaces.
\end{remark}

The next theorem follows from two cited results and gives an estimate for the embeddable genus of a neighborhood of a regular multibranched surface.

\begin{theorem}[{\cite[Theorem 3.6]{MO}}, {\cite[Lemma 2.2]{EO2}}]
For a regular multibranched surface $X$ and any neighborhood $N \in \mathcal{N}(X)$, the next inequality holds.
\[
rank H_1(X) -g(\partial N)\le eg(N) \le rank H_1(G_N) + g(\partial N),
\]
where $G_N$ denotes the abstract dual graph of $N$ and $g(\partial N)$ denotes the sum of genera of all components of $\partial N$.
\end{theorem}

\subsection{Lower bounds}

The following lower bounds for the minimum genus and the maximum genus are known.

\begin{theorem}[\cite{EO2}, cf. {\cite[Theorem 1.3]{T}}]\label{lower}
For a regular multibranched surface $X$, the next inequalities hold.
\begin{eqnarray}
\min g(X) &\ge& rank H_1(X) - \max_{N \in \mathcal{N}(X)} g(\partial N)\\
\max g(X) &\ge& rank H_1(X) - \min_{N \in \mathcal{N}(X)} g(\partial N)
\end{eqnarray}
\end{theorem}

\subsection{Graph product $G \times S^1$}

For a graph $G$, we obtain a regular multibranched surface by taking a product with $S^1$.
We consider the genus of a regular multibranched surface which forms $G\times S^1$, and by using Theorem \ref{lower}, we obtain the following theorem which is an interplay of the genus of a graph $G$ and the genus of a multibranched surface $G\times S^1$.

The {\em minimum genus} $\min g(G)$ of a graph $G$ is defined as the minimal genus of closed orientable surfaces in which $G$ can be embedded.
The {\em maximum genus} $\max g(G)$ of a graph $G$ is defined as the maximal genus of closed orientable surfaces in which $G$ can be embedded and the complement of $G$ consists of open disks.
It is remarkable that Xuong and Nebesk\'{y} determined the maximum genus of a graph by a completely combinatorial formula ({\cite[Theorem 3]{X}}, {\cite[Theorem 2]{N81}}).

\begin{theorem}[{\cite[Corollary 1.2]{T}}, \cite{EO2}]\label{product}
For a graph $G$, the next equalities hold.
\begin{eqnarray}
\min g(G\times S^1) &=& 2 \min g(G)\\
\max g(G\times S^1) &=& 2 \max g(G)
\end{eqnarray}
\end{theorem}

In {\cite[Corollary 1.2]{T}}, it was shown that the minimal number $\dim H_1(M ;F)$ for a closed orientable 3-manifolds $M$ containing $G\times S^1$ equals to $2\min g(G)$, where $F=\Bbb{Z}_p$ or $\Bbb{Q}$.
It is well-known that $g(M)\ge \dim H_1(M ;F)$.
Hence the inequality $\min g(G\times S^1) \ge 2 \min g(G)$ in Theorem 3.5 (3.3) holds.

\subsection{Spine of closed 3-manifolds}

A multibranched surface $X$ is called a {\em 2-stratifold} if every prebranch $c$ of $X$ satisfies $d(c)>2$.
G\'{o}mez-Larra\~{n}aga--Gonz\'{a}lez-Acu\~{n}a--Heil studied 2-stratifolds from a view point of 3-manifold groups.
They asked the following questions.

\begin{question}
Which 3-manifolds $M$ have fundamental groups isomorphic to the fundamental group of a 2-stratifold?
\end{question}

\begin{question}
Which closed 3-manifolds $M$ have spines that are 2-stratifolds? 
\end{question}

Recall that a subpolyhedron $P$ of a 3-manifold $M$ is a spine of $M$, if $M-{Int}(B^3)$ collapses to $P$, where $B^3$ is a 3-ball in $M$. An equivalent definition is that $M-P$ is homeomorphic to an open 3-ball.
They completely answered these questions.

\begin{theorem}[{\cite[Theorem 1]{GGH}}]
Let $M$ be a closed 3-manifold and $X_G$ be a 2-stratifold. If $\pi_1(M)\cong \pi_1(X_G)$, then $\pi_1(M)$ is a free product of groups, where each factor is cyclic or $\mathbb{Z}\times \mathbb{Z}_2$.
\end{theorem}

\begin{theorem}[{\cite[Theorem 2]{GGH}}]
A closed 3-manifold $M$ has a 2-stratifold as a spine if and only if $M$ is a connected sum of lens spaces, $S^2$-bundles over $S^1$, and $P^2\times S^1$'s.
\end{theorem}

\subsection{Neighborhood equivalence}

In this subsection, we assume that a multibranched surface is regular, does not have disk sectors, and the degree is greater than 2 for each branch.
Let $A$ be an annulus sector of $X$ whose boundary consists of two branches, where at least one branch has the wrapping number $1$, or a M\"{o}bius-band sector whose boundary has the wrapping number $1$.
An {\em IX-move} along $A$ is an operation shrinking $A$ into the core circle, and an {\em XI-move} is a reverse operation of an IX-move.

If two multibranched surfaces $X,X'$ embedded in a 3-manifold $M$ are related by a IX-moves or XI-moves, then the regular neighborhoods $N(X), N(X')$ are isotopic in $M$.
The following theorem states that the converse holds.

\begin{theorem}[\cite{IKOS}]\label{IH}
Let $X,\ X'$ be multibranched surfaces embedded in an orientable $3$-manifold $M$.
If $N(X)$ is isotopic to $N(X')$ in $M$, then $X$ is transformed into $X'$ by a finite sequence of IX-moves, XI-moves and isotopies.
\end{theorem}

For more broad class, Matveev--Piergallini theorem is known:
Two simple polyhedra embedded in a 3-manifold has isotopic neighborhoods if and only if they are connected by a sequence of {\it $2 \leftrightarrow 3$ moves}, {\it $0 \leftrightarrow 2$ moves} moves and isotopy (\cite{M88}, \cite{P88}).

\subsection{Neighborhood partial order}\label{NE}

Let $X$ be an essential multibranched surface embedded in a closed orientable 3-manifold $M$.
We say that a sector $s$ is {\em excess} if it is boundary-parallel in $M-int N(X-s)$.
A multibranched surface $X$ is said to be {\em efficient} if every sector is not excess.

In this subsection, we restrict multibranched surfaces to the set $\mathcal{X}$ of all connected compact multibranched surfaces $X$ embedded in a closed orientable 3-manifold $M$ satisfying the following conditions:
$X$ is maximally spread (that is, applied XI-moves to $X$ as much as possible), essential and efficient in $M$, has no open disk sector, no branch of degree 1 or 2.

Under the influence of Theorem \ref{IH}, we define an equivalence relation on $\mathcal{X}$ as follows.
Two multibranched surfaces $X$ and $X'$ in $\mathcal{X}$ are {\em neighborhood equivalent}, denoted by $X\overset{\mathrm{N}}\sim X'$, if $X$ is transformed into $X'$ by a finite sequence of IH-moves.
Moreover, we define a binary relation $\le$ over $\mathcal{X}$ as follows.

\begin{definition}\label{relation}
For $X,\, Y\in \mathcal{X}$, we denote $X \le Y$ if 
\begin{enumerate}
  \setlength{\itemsep}{0cm}
\item there exists an isotopy of $Y$ in $M$ so that $Y\subset N(X)$ and $B_Y\subset N(B_X)$, and
\item there exists no essential annulus in $N(X)-Y$.
\end{enumerate}
\end{definition}

For equivalence classes $[X], [Y]\in \mathcal{X}/\overset{\mathrm{N}}\sim $, we define a {\em neighborhood partial order} $\preceq$ over $\mathcal{X}/\overset{\mathrm{N}}\sim $ so that $[X]\preceq [Y]$ if $X\le Y$.

\begin{theorem}[\cite{O}]\label{partial}
The relation $\preceq$ is well-defined on $\mathcal{X}/\overset{\mathrm{N}}\sim $ and 
$(\mathcal{X}/\overset{\mathrm{N}}\sim ; \preceq)$ is a partially ordered set.
\end{theorem}

We say that $B_X$ is {\em toroidal} if there exists an essential torus $T$ in the exterior $E(B_X)$ of $B_X$ in $M$, that is, $T$ is incompressible in $E(B_X)$ and $T$ is not parallel to a torus in $\partial E(B_X)$.
We say that $E_X$ is {\em cylindrical} if there exists an essential annulus $A$ with $A\cap X =A\cap E_X=\partial A$, that is, $A$ is incompressible and $A$ is parallel to neither an annulus in $E_X$ nor an annulus in $\partial E(B_X)$.

\begin{theorem}[\cite{O}]\label{sufficient}
For equivalence classes $[X], [Y]\in \mathcal{X}/\overset{\mathrm{N}}\sim $, 
if $[X]\preceq [Y]$ and $[X]\ne [Y]$, then either $B_Y$ is toroidal or $S_Y$ is cylindrical.
\end{theorem}

Theorem \ref{sufficient} provides a sufficient condition for an equivalent class $[X]\in \mathcal{X}/\overset{\mathrm{N}}\sim $ to be minimal with respect to the partial order of $(\mathcal{X}/\overset{\mathrm{N}}\sim ; \preceq)$, that is, if $B_X$ is atoroidal and $E_X$ is acylindrical, then $[X]$ is minimal.


\subsection{Essential decomposition | Eudave-Mu\~{n}oz knots type}

Let $X$ be a multibranched surface embedded in the 3-sphere $S^3$, which decomposes $S^3$ into regions $V_1,\ldots, V_n$.
If $X$ is essential, then we call this decomposition $S^3=V_1\cup\cdots\cup V_n$ an {\em essential decomposition}.
As explained in Subsection \ref{essential}, a link with an essential surface of non-meridional boundary slope gives an essential decomposition.

In this subsection, we recall Eudave-Mu\~{n}oz knots (\cite{EM}) in the language of multibranched surface.
Let $X$ be a multibranched surface which has a twice punctured torus as a sector $s$ and a single branch $l$, where one prebranch $c$ has $od(c)=2$ and another prebranch $c'$ has $od(c')=-2$.

Suppose that $X$ is embedded in $S^3$ so that it is essential and two regions of $S^3-X$ are genus two handlebodies, say $H$ and $W$.
Then, by combining \cite{EM} with \cite{GL}, the branch $l$ forms an Eudave-Mu\~{n}oz knot.
From the point of view that any essential embedding restricts the knot type of the branch,
this phenomenon is especial in low-dimensional geometric topology.

\begin{figure}[htbp]
\begin{center}
\includegraphics[width=.25\linewidth]{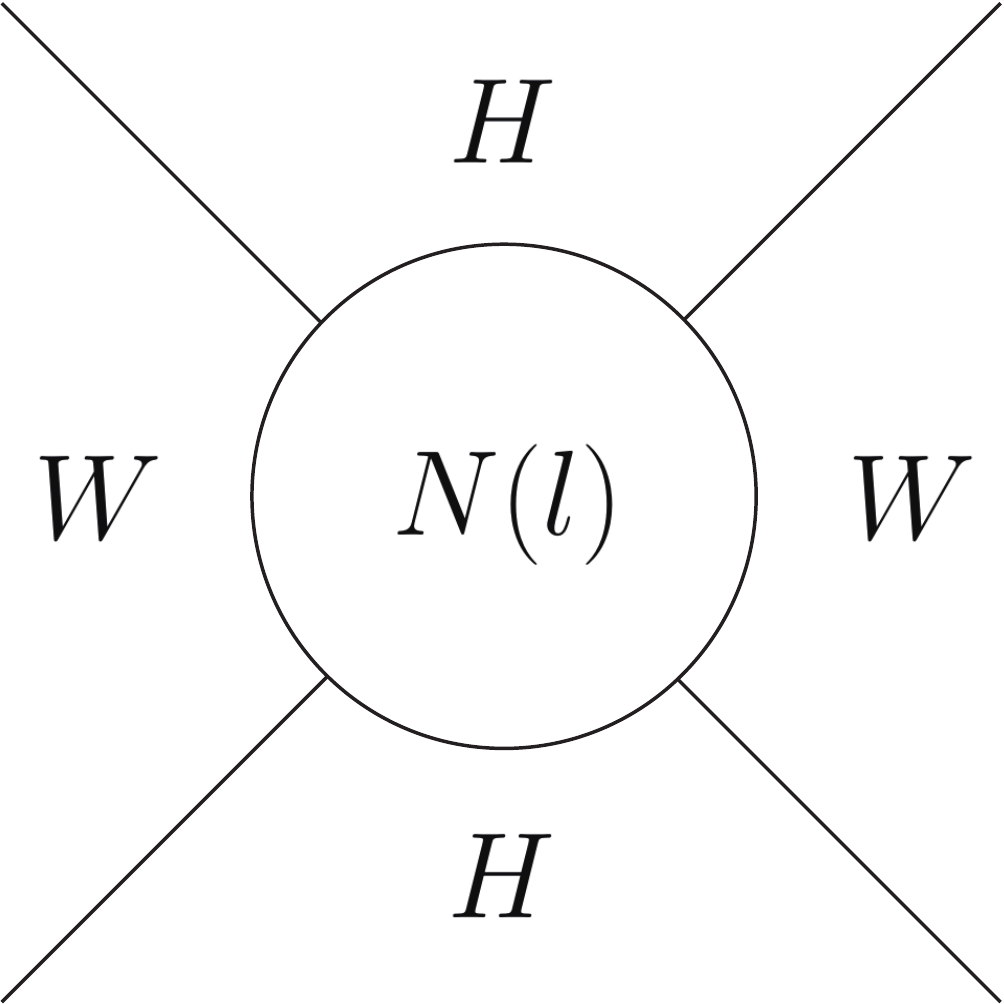}
\end{center}
	\caption{$(1,2,2;2)$-trisection coming from Eudave-Mu\~{n}oz knots}
	\label{EM}
\end{figure}

Eudave-Mu\~{n}oz knots appears in the last piece of the classification of essential annuli in the exterior of genus two handlebody-knots in $S^3$ (\cite{KO2}).
We take a regular neighborhood $N(l)$ and denote two handlebodies $S^3-N(l)-s$ by $H$ and $W$ again.
See Figure \ref{EM} for the configulation.
Put $A=N(l)\cap W$.
Then, $H$ is a genus two handlebody-knots with an essential annulus $A$ of Type 4 in \cite{KO2}.

This configulation as in Figure \ref{EM} also provides a nice example of trisection.
Let $X'$ be a multibranched surface which has two branches $b\cup b'=N(l)\cap H\cap W$ and three sectors $s_1=H\cap W$, $s_2=N(l)\cap H$ and $s_3=N(l)\cap W$.
Then, $X'$ gives an essential decomposition $S^3=N(l)\cup H\cup W$, where the triple of genera of three handlebodies is $(1,2,2)$ and the number of branches is $2$.
Thus this gives a $(1,2,2;2)$-trisection of $S^3$.
Moreover, it is shown in {\cite[Proposition 4.7.1]{Ko}} that this trisection is not stabilizations of any other trisection.

\subsection{Efficient embedding | Universal bounds}

Recall the relation between essential surfaces in link exteriors and essential multibranched surfaces in Subsection \ref{essential}, and  the definition of efficient embedding in Subsection \ref{NE}.
Suppose that $X$ is an essential and efficient multibranched surface embedded in a 3-manifold.
Then, we have a link and essential surfaces in the link exterior, and moreover any pair of essential surfaces is not mutually parallel.

Let $X$ be a multibranched surface which has a single branch and $n$ once punctured tori each of which has an oriented degree 1.
Suppose that $X$ is embedded in $S^3$ so that it is essential and efficient and the branch forms a hyperbolic knot.
Then, we have a hyperbolic knot bounding $n$ genus 1 Seifert surfaces which are not mutually parallel.  
Tsutsumi first showed that the number $n$ is at most 7 (\cite{Tsu}).
After that, Eudave-Mu\~{n}oz -- Ram\'{i}rez-Losada -- Valdez-S\'{a}nchez showed that $n$ is at most 6 and provided an example of such embedding of $X$ for $n=5$ (\cite{ERV}).
Finally, Valdez-S\'{a}nchez showed that $n$ is at most 5 (\cite{V}) and therefore this bound is optimal.


This phenomenon is also especial in low-dimensional geometric topology.
In general, contrary to the above, there is no upper bound.
Tsutsumi showed that for any positive integer $n$, there is a genus one hyperbolic knot in $S^3$ which bounds mutually non-parallel incompressible Seifert surfaces $S, F_1,\ldots, F_n$ where $S$ is genus one and $F_i$ is genus two (\cite[Theorem 5.5]{Tsu}).

\section{Forbidden minors for $S^3$}\label{forbidden}

\subsection{Minor}

In this subsection, we allow the degree $deg(B_i)$ of a branch $B_i$ to be $1$ or $2$.

We denote by $\mathcal{M}$ the set of all regular multibranched surfaces (modulo homeomorphism).
For $X$, $Y\in \mathcal{M}$,
we write $X < Y$ if $X$  is obtained by removing a sector of $Y$, or $X$ is obtained from $Y$ by an IX-move.  
We define an equivalence relation $\sim$ on $\mathcal{M}$ as follows: if $X < Y$ and $Y < X$, then $X \sim Y$.

We define a partial order $\prec$ on $\mathcal{M}/\sim$ as follows. 
Let $X$, $Y \in \mathcal{M}$. 
We denote $[X] \prec [Y]$ if there exists a finite sequence $X_1, \ldots, X_n \in \mathcal{M}$ such that $X_1 \sim X$, $X_n \sim Y$ and $X_1 < \cdots < X_n$.

\subsection{Obstruction set}

A multibranched surface class $[X]$ is called a {\em minor} of another multibranched surface class $[Y]$ if $[X] \prec [Y]$. 
In particular, $[X]$ is called a {\em proper minor} of $[Y]$ if $[X] \prec [Y]$ and $[Y] \not= [X]$.
A subset $\mathcal{P}$ of $\mathcal{M}/\sim$ is said to be {\em minor closed} if for any {$[X] \in \mathcal{P}$, every minor of $[X]$ belongs to $\mathcal{P}$}.
For a minor closed set $\mathcal{P}$, we define the {\em obstruction set} $\Omega(\mathcal{P})$ by all elements $[X] \in \mathcal{M}/\sim $ such that $[X] \not \in \mathcal{P}$ and every proper minor of $[X]$ belongs to $\mathcal{P}$.

The set of multibranched surfaces embeddable into $S^3$, denoted by $\mathcal{P}_{S^3}$, is minor closed.
As a 2-dimensional version of Kuratowski's and Wagner's theorems,
we consider the next problem.

\begin{problem}\label{obstruction_set}
Characterize the obstruction set $\Omega(\mathcal{P}_{S^3})$.
\end{problem}

We summarize all known results on $\Omega(\mathcal{P}_{S^3})$ at the present moment.
As we shall see later, (2) and (3) in Theorem \ref{obstruction} are infinite family of multibranched surfaces.

\begin{theorem}\label{obstruction}
The following multibranched surfaces belong to $\Omega(\mathcal{P}_{S^3})$.
\begin{enumerate}
  \setlength{\itemsep}{0cm}
\item $K_5\times S^1$ and $K_{3,3}\times S^1$ $($\cite{S}$)$
\item $X_1$, $X_2$, $X_3$ $($\cite{EMO}$)$
\item $X_g(p_1,\ldots,p_n)$ $($\cite{MO}$)$
\end{enumerate}
\end{theorem}


\begin{remark}
(1) Since any proper minor of $K_5$ and $K_{3,3}$ is planar, any proper minor of $K_5\times S^1$ and $K_{3,3}\times S^1$ can be embedded in $D^2\times S^1\subset S^3$.

(2) We say that a multibranched surface $X$ is {\em critical} for $S^3$ if $X$ cannot be embedded in $S^3$ and for any $x\in X$, $X-x$ can be embedded in $S^3$.
It is shown in \cite{EMO} that $X_1$, $X_2$, $X_3$ are critical for $S^3$.

(3) Since $X_g(p_1,\ldots,p_n)$ has a single sector, the minimality for $\mathcal{P}_{S^3}$ naturally holds.
\end{remark}

Theorem \ref{obstruction} (1) was proved in \cite[Theorem 1]{S}.
It also follows Theorem \ref{product} and Kuratowski's and Wagner's theorem.

The families of multibranched surfaces $X_1,\ X_2,\ X_3$ in Theorem \ref{obstruction} (2) are given as follows.

\begin{figure}[htbp]
	\begin{center}
	\includegraphics[width=.2\linewidth]{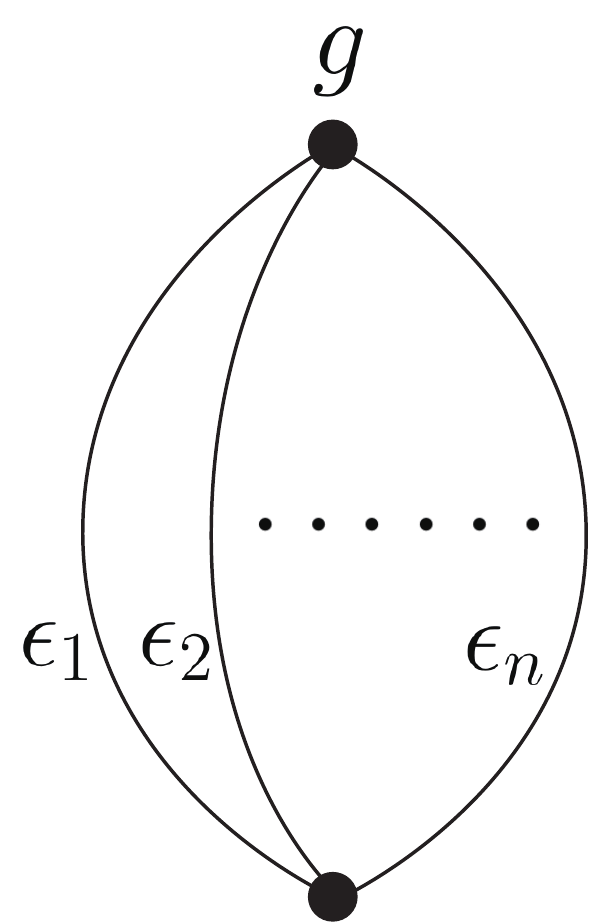}
	\end{center}
	\caption{A graph representation of $X_1$}
	\label{X_1}
\end{figure}

Let $X_1$ be a multibranched surface which is obtained from a single sector of genus $g$, $n$ boundary components and a single branch by a covering map with degree $\epsilon_i$ on each prebranch.
See Figure \ref{X_1} for a graph representation.
We assume that $\epsilon_i=\pm p$ for the regularity of $X_1$.
Then the incident matrix is $M_{X_1}=\big(\sum_{i=1}^n \epsilon_i\big)$.
If $\big|\sum_{i=1}^n \epsilon_i\big|>1$, then $H_1(X_1)$ has a torsion and $X_1$ cannot be embedded in $S^3$.
Conversely, if $\big|\sum_{i=1}^n \epsilon_i\big|\le 1$, then by \cite[Theorem 3.2]{EMO}, $X_1$ can be embedded in $S^3$.
Hence $X_1\in \Omega(\mathcal{P}_{S^3})$ if and only if $\big|\sum_{i=1}^n \epsilon_i\big|>1$.

\begin{figure}[htbp]
	\begin{center}
	\includegraphics[width=.4\linewidth]{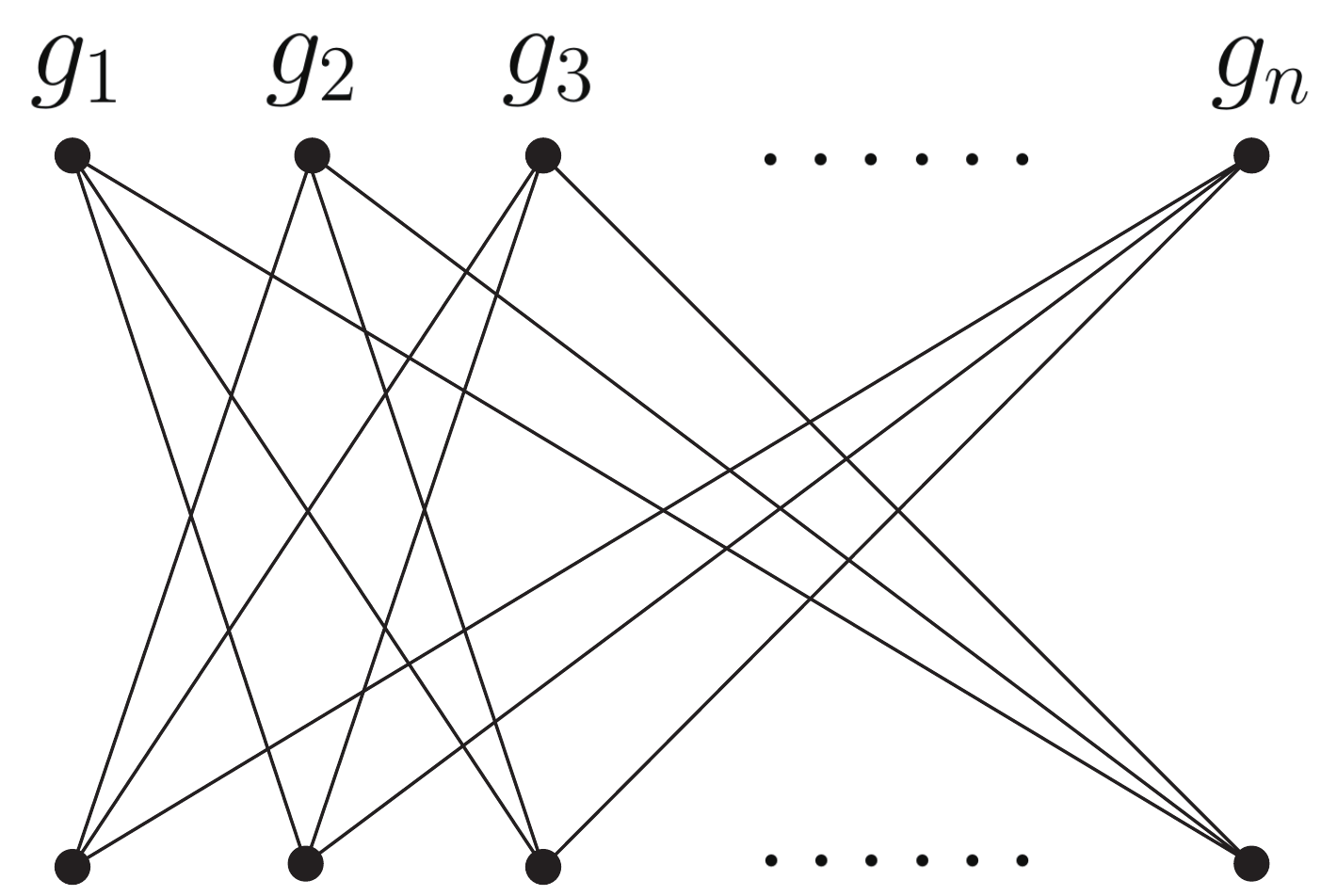}
	\end{center}
	\caption{A graph representation of $X_2$}
	\label{X_2}
\end{figure}

Let $X_2$ be a multibranched surface which has a graph representation as shown in Figure \ref{X_2}, where $n\ge 3$, all wrapping number is 1 (we omit the labels on edges).
Then by \cite[Theorem 3.3]{EMO}, $X_2\in \Omega(\mathcal{P}_{S^3})$.

The incident matrix of $X_2$ is:
\[
M_{X_2}=\left( \begin{array}{cccccc}
0 & 1 & \cdots & \cdots & 1 \\
1 & 0 & \ddots & \ddots & \vdots \\
\vdots & \ddots & \ddots & \ddots & \vdots \\
\vdots & \ddots & \ddots & \ddots & 1 \\
1 & \cdots & \cdots & 1 & 0
\end{array} \right)
\]
Since $det(M_{X_2})=(-1)^{n+1}(n-1)$ and $n\ge 3$, $X_2$ has a torsion.

\begin{figure}[htbp]
\begin{center}
\includegraphics[width=.45\linewidth]{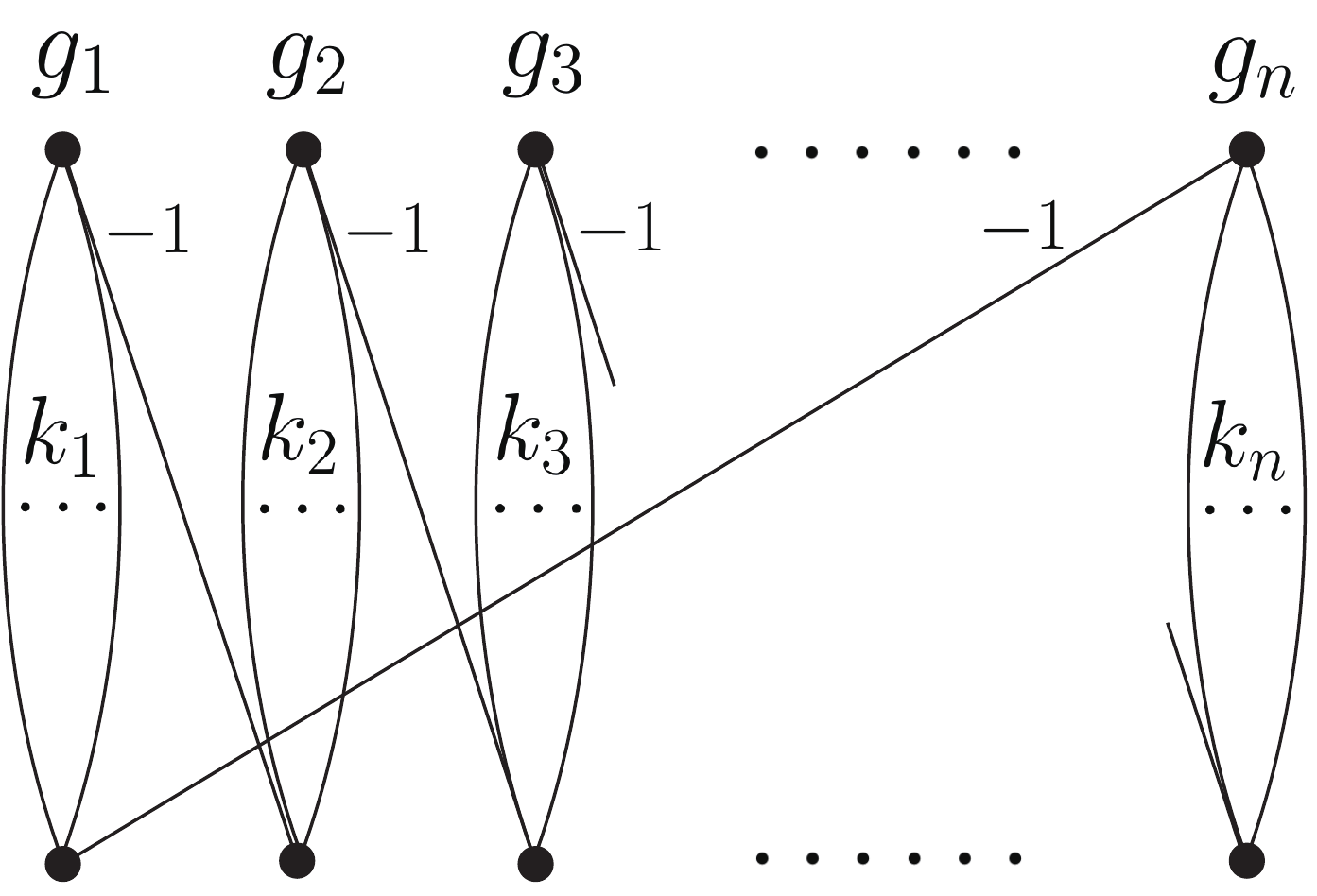}
\end{center}
	\caption{A graph representation of $X_3$}
	\label{X_3}
\end{figure}

Let $X_3$ be a multibranched surface which has a graph representation as shown in Figure \ref{X_3}, where $n\ge 2,\ k_i\ge 1,\ k_1k_2k_3\cdots k_n\ge3$, all wrapping number is 1 unless otherwise specified.
Then by \cite[Theorem 3.7]{EMO}, $X_3\in \Omega(\mathcal{P}_{S^3})$.

The incident matrix of $X_3$ is:
\[
M_{X_3}=\left( \begin{array}{cccccc}
k_{1} & 0 & \cdots & \cdots & 0 & -1 \\
-1 & k_{2} & 0 & \ddots & \ddots & 0\\
0 & -1 & k_3 & \ddots & \ddots & \vdots\\
\vdots & 0 & -1 & \ddots & \ddots & 0\\
\vdots & \ddots & \ddots & \ddots & k_{n-1} & 0\\
0 & \cdots & \cdots & 0 & -1 & k_n
\end{array} \right)
\]
Since $det(M_{X_3})=k_1k_2k_3\cdots k_n -1 \geq 2$,  $X_3$ has a torsion.

\begin{figure}[htbp]
\begin{center}
\includegraphics[width=.3\linewidth]{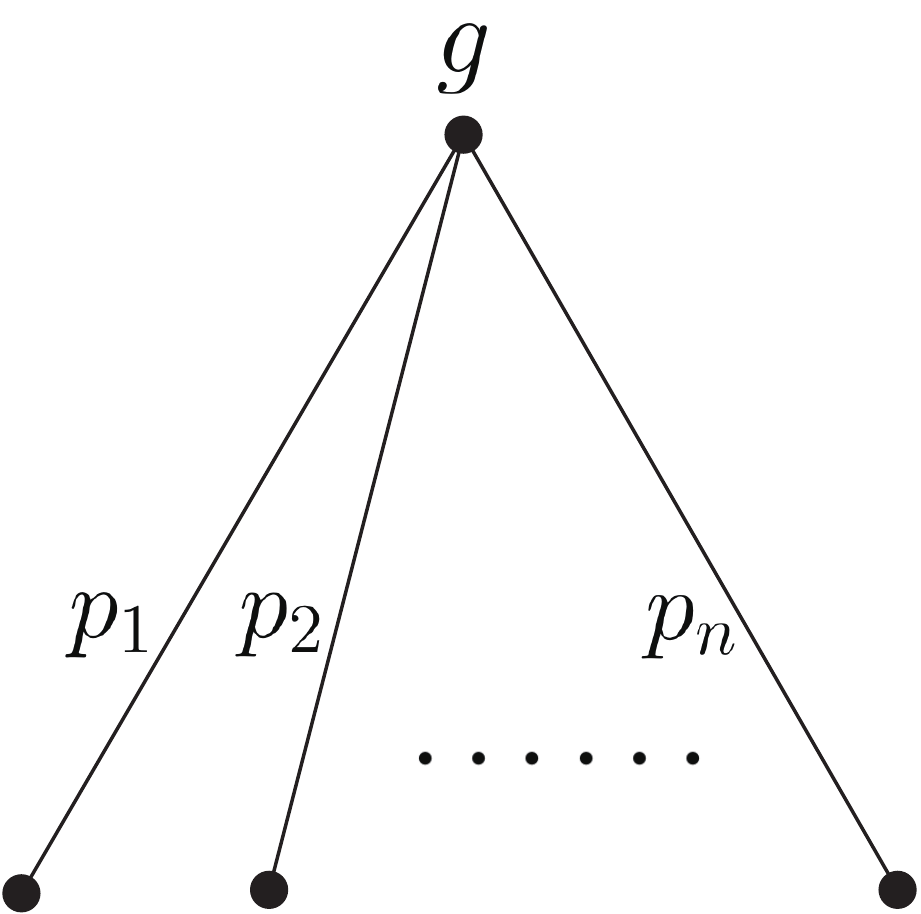}
\end{center}
	\caption{A graph representation of $X_g(p_1,\ldots,p_n)$}
	\label{X_g}
\end{figure}

The multibranched surface $X_g(p_1,\ldots,p_n)$ in Theorem \ref{obstruction} (3) was first presented in \cite[Example 4.3]{MO}.
Let $X_g(p_1,\ldots,p_n)$ be a multibranched surface which has a graph representation as shown in Figure \ref{X_g}, where $n\ge 1$, $p=gcd\{ p_1, \ldots, p_n\}>1$.
As we have seen in Example \ref{crosscap}, a closed non-orientable surface of crosscap number $n$ is homeomorphic to  $X_0(2,\ldots,2)$.

As shown in \cite[Example 4.3]{MO}, we have $H_1(X_g(p_1,\ldots,p_n))=\left( \mathbb{Z}/p \mathbb{Z} \right) \oplus \mathbb{Z}^{2g+n-1}$.
Hence, $X_g(p_1,\ldots,p_n)$ cannot be embedded in $S^3$ since $p>1$.

\subsection{Beyond torsion}\label{beyond}

In the previous subsection, we conclude that some multibranched surfaces cannot be embedded in $S^3$ in virtue of torsion part of the first homology group.
Then, the following inverse problem naturally arises.

\begin{problem}[{\cite[Problem, p.631]{MO}}]\label{inverse}
If $p=1$, then can $X_g(p_1,\ldots,p_n)$ be embedded in $S^3$?
\end{problem}

The following theorem gives a partial answer to Problem \ref{inverse}.

\begin{theorem}[{\cite[Theorem 1.5]{EO}}]
If $p=1$, then $X_g(p_1, p_2, p_3)$ can be embedded in $S^3$ for a sufficiently large $g$.
\end{theorem}

But, how can we say about Problem \ref{inverse} when $g=0$?
This is related with a main thema in \cite{EO}.
In \cite{EO}, we characterized non-hyperbolic 3-component links in the 3-sphere whose exteriors contain essential 3-punctured spheres with non-integral boundary slopes.
This implies that we can derive a formula for the triple $p_1, p_2, p_3$ ({\cite[Proposition 1.4]{EO}}).
For hyperbolic links, we conjectured the following.

\begin{conjecture}[{\cite[Conjecture 1.1]{EO}}, c.f. \cite{GS1986}, \cite{GL1987}]\label{cabling}
There does not exist an essential n-punctured sphere with non-meridional, non-integral boundary slope in a hyperbolic link exterior in the 3-sphere.
\end{conjecture}

It can be checked that the triple $(5, 7, 18)$ does not satisfy the formula in {\cite[Proposition 1.4]{EO}}.
Therefore, assuming Conjecture \ref{cabling}, we conclude that $X_0(5, 7, 18)$ cannot be embedded in $S^3$.

On the other hand, if we allow embeddings in 3-manifolds other than $S^3$, then Problem \ref{inverse} holds.
We use a result in \cite{GA} that a compact $3$-manifold $M$ with connected boundary can be embedded in a homology $3$-sphere if and only if $H_1(M)$ is free and $H_2(M)=0$.
Since for a unique neighborhood $N \in \mathcal{N}(X_g(p_1,\ldots,p_n))$, $H_1(N)$ is free and $H_2(N)=0$ when $p=1$, we have the following.

\begin{theorem}
If $p=1$, then $X_g(p_1,\ldots,p_n)$ can be embedded in a homology $3$-sphere.
\end{theorem}







\section{The prospects for multibranched surfaces}

The author would like to close this survey article by stating the following prospects.

Firstly, it is important to characterize essential and efficient decompositions of $S^3$, where we say that a decomposition $S^3=V_1\cup \cdots \cup V_n$ by a multibranched surface $X$ is {\em efficient} if $X$ is efficient.
It can be applied to Polycontinuous patterns, Trisections, Essential surfaces as stated in Introduction.

Secondly, it is a fundamental problem to characterize the obstruction set $\Omega(\mathcal{P}_{S^3})$.
This problem has a difficulty as stated in Subsection \ref{beyond}, but it has also an interest in Conjecture \ref{cabling}.

\bibliographystyle{amsplain}

\bigskip
\noindent{\bf Acknowledgements.}
The author would like to thank to Fico Gonz\'{a}lez-Acu\~{n}a and Arkadiy Skopenkov for informimg me related results.

\end{document}